\documentclass[oneside,reqno]{amsart}

\usepackage[active]{srcltx}

 \makeatletter
\renewcommand*\subjclass[2][2010]{%
  \def\@subjclass{#2}%
  \@ifundefined{subjclassname@#1}{%
    \ClassWarning{\@classname}{Unknown edition (#1) of Mathematics
      Subject Classification; using '2010'.}%
  }{%
    \@xp\let\@xp\subjclassname\csname subjclassname@#1\endcsname
  }%
}

 \makeatother
\newtheorem{theorem}{Theorem}[section]

\theoremstyle{definition}

\numberwithin{equation}{section}
 \makeatletter
\renewcommand*\subjclass[2][2010]{%
  \def\@subjclass{#2}%
  \@ifundefined{subjclassname@#1}{%
    \ClassWarning{\@classname}{Unknown edition (#1) of Mathematics
      Subject Classification; using '1991'.}%
  }{%
    \@xp\let\@xp\subjclassname\csname subjclassname@#1\endcsname
  }%
}
 \makeatother
\begin{document}

\title[Proof of  a  congruence for harmonic numbers\dots]{Proof of  
a  congruence for harmonic numbers conjectured by Z.-W. Sun}

\author{Romeo Me\v strovi\' c}
\address{Department of Mathematics,
Maritime Faculty, University of Montenegro, 
Dobrota 36, 85330 Kotor, Montenegro 
romeo@ac.me}

\begin{abstract}
 For a positive integer $n$ let $H_n=\sum_{k=1}^{n}1/k$
be the $n$th  harmonic number.
In this note we prove that for any prime $p\ge 7$,
   $$
\sum_{k=1}^{p-1}\frac{H_k^2}{k^2}
\equiv\frac{4}{5}pB_{p-5}\pmod{p^2},
   $$
which confirms the conjecture recently proposed by Z.-W. Sun. 
Furthermore, we also prove two similar congruences modulo
$p^2$.  
   \end{abstract}

  \maketitle

\section{{\bf The Main Result and Its Proof}}

For a positive integer $n$, 
the $n$th  {\it harmonic number} $H_n$
is defined as
 $$
H_n=\sum_{k=1}^n\frac{1}{k}
  $$   
(in addition, we define  $H_0=0$).

Recently, Zhi-Wei  Sun \cite{s}  obtained  basic congruences modulo a prime
$p\ge 5$ for several sums of terms involving harmonic numbers.
In particular,  Sun established
$\sum_{k=1}^{p-1}H_k^r\,(\bmod{\, p^{4-r}})$ for $r=1,2,3$. 
Further generalizations of these congruences have been recently obtained by 
Tauraso in \cite{t}. Notice also that Sun \cite[Theorem 1.1]{s}
proved that for any prime $p>5$
  \begin{equation}\label{con1}
\sum_{k=1}^{p-1}\frac{H_k^2}{k^2}\equiv 0\pmod{p},
   \end{equation}
and conjectured \cite[the second part of Conjecture 1.1]{s} 
(also see \cite[Conjecture A37 (v), p. 31]{s2}) that for such a prime $p$,
  \begin{equation}\label{con2}
\sum_{k=1}^{p-1}\frac{H_k^2}{k^2}
\equiv\frac{4}{5}pB_{p-5}\pmod{p^2},
   \end{equation}
where $B_0,B_1,B_2,\ldots,$ are Bernoulli numbers given by
   $$
B_0=1\quad\mathrm{and}\quad \sum_{k=0}^{n}{n+1\choose k}B_k=0\,\, 
(n=1,2,3,\ldots).
  $$
Recall also that the first part of Conjecture 1.1 in \cite{s} 
was recently confirmed by Z.-W. Sun and L. L. Zhao in
\cite[the congruence (\ref{con1}) of Theorem 1.1]{sz}; for another proof see
\cite[the congruence (7) of Corollary 1.2]{m}.

In this note we prove the congruence (\ref{con2}) and 
two similar congruences contained in the following result.
\begin{theorem}  Let $p>5$ be a prime. Then
  \begin{equation}\label{con3}
\sum_{k=1}^{p-1}\frac{H_k^2}{k^2}
\equiv\frac{4}{5}pB_{p-5}\pmod{p^2},
   \end{equation}
\begin{equation}\label{con4}
\sum_{k=1}^{p-1}\frac{H_k^3}{k}
\equiv\frac{3}{2}pB_{p-5}\pmod{p^2},
   \end{equation}
\begin{equation}\label{con5}
\sum_{k=1}^{p-1}\frac{H_k}{k^3}
\equiv -\frac{1}{10}pB_{p-5}\pmod{p^2}.
   \end{equation}
\end{theorem}  

The proof of Theorem 1.1 presented below is 
based on some  congruences for multiple harmonic sums obtained by 
J. Zhao \cite{z}.    

\begin{proof}[Proof of Theorem {\rm 1.1.}]

For simplicity, here we denote
   $$
H(2,1,1)=\sum_{1\le i<j< k\le p-1}\frac{1}{i^2jk},\quad
H(1,2,1)=\sum_{1\le i<j< k\le p-1}\frac{1}{ij^2k},
    $$
  $$
H(m,n)=\sum_{1\le i<j\le  p-1}\frac{1}{i^mj^n}, 
\,\,\mathrm{for}\,\, m,n=1,2,\ldots,\,\,
H(n)=\sum_{k=1}^{p-1}\frac{1}{k^n}\,\,\mathrm{for}\,\, n=1,2,\ldots.
  $$

{\it Proof of {\rm (\ref{con3})}}.
By  {\it Wolstenholme theorem} (see e.g., \cite[Theorem 1]{a} or 
\cite[p. 255]{gr}), for any prime $p>3$,  
$H(1):=\sum_{k=1}^{p-1}1/k\equiv 0\,(\bmod{\,p^2})$,
or equivalently, for each $j=1,2,\ldots,p-2$ the following congruence 
holds
 \begin{equation}\label{con6}
\frac{1}{j+1}+\frac{1}{j+2}+\cdots +\frac{1}{p-1}\equiv -\left(
1+\frac{1}{2}+\cdots +\frac{1}{j-1}+\frac{1}{j}\right)\pmod{p^2}.
 \end{equation}
Applying  the congruence (\ref{con6}), we find that 
    \begin{eqnarray}
&&H(1,2,1)=\sum_{1\le i<j< k\le p-1}\frac{1}{ij^2k}
=\sum_{j=2}^{p-1}\frac{1}{j^2}\sum_{i=1}^{j-1}\frac{1}{i}
\sum_{k=j+1}^{p-1}\frac{1}{k}\nonumber\\
&=&\sum_{j=2}^{p-1}\frac{1}{j^2}
\left(1+\frac{1}{2}+\cdots +\frac{1}{j-1}\right)
\left(\frac{1}{j+1}+\frac{1}{j+2}+\cdots +\frac{1}{p-1}\right)\nonumber\\
&\equiv& \sum_{j=2}^{p-1}\frac{1}{j^2}
\left(1+\frac{1}{2}+\cdots +\frac{1}{j-1}\right)
\left(-\left(1+\frac{1}{2}+\cdots +\frac{1}{j}\right)\right)\pmod{p^2}
\nonumber\\
&=&\sum_{j=1}^{p-1}\frac{1}{j^2}\left(H_{j}-\frac{1}{j}\right)
\left(-H_j\right)=-\sum_{j=1}^{p-1}\frac{H_j^2}{j^2}+
\sum_{j=1}^{p-1}\frac{H_j}{j^3}\nonumber\\
&=&-\sum_{j=1}^{p-1}\frac{H_j^2}{j^2}+\sum_{j=1}^{p-1}\frac{H_{j-1}
+\frac{1}{j}}{j^3}=
-\sum_{j=1}^{p-1}\frac{H_j^2}{j^2}+\sum_{1\le i<j\le p-1}\frac{1}{ij^3}
+\sum_{j=1}^{p-1}\frac{1}{j^4}\nonumber\\
&=&-\sum_{j=1}^{p-1}\frac{H_j^2}{j^2}+H(1,3)+H(4)\hfill
\qquad\qquad\qquad\qquad\qquad\quad\pmod{p^2}.\label{con7}
   \end{eqnarray}

From the shuffle relation  (see e.g., \cite[(3.6) on page 89]{z})
  $$
H(1)H(3)=H(1,3)+H(3,1)+H(4)
  $$
and the fact that $p^2\mid H(1)$  it follows that
  $$
H(1,3)+H(4)\equiv -H(3,1)\pmod{p^2}.
  $$
Substituting this into the right hand side of (\ref{con7}), 
we immediately obtain
  \begin{equation}\label{con8}
H(1,2,1)+H(3,1)\equiv -\sum_{j=1}^{n}\frac{H_j^2}{j^2}
 \pmod{p^2}.
  \end{equation} 
Further, by the shuffle relation  
  $$
H(1)H(2,1)=2H(2,1,1)+H(1,2,1)+H(3,1)+H(2,2)
  $$
again using the fact that $p^2\mid H(1)$, we get
        \begin{equation}\label{con9}
    H(1,2,1)+H(3,1)\equiv -2H(2,1,1)-H(2,2)\pmod{p^2}
   \end{equation}
By the congruence (3.22) in \cite[Proposition 2.7, p. 93]{z}, we have
\begin{equation}\label{con10}
H(2,1,1)\equiv \frac{3}{5}pB_{p-5}\pmod{p^2},
   \end{equation}
and by \cite[Theorem 3.2, p. 88]{z} or \cite{zc},
\begin{equation}\label{con11}
H(2,2)\equiv -\frac{2}{5}pB_{p-5}\pmod{p^2}.
   \end{equation}
From (\ref{con9}), (\ref{con10}) and (\ref{con11}) we find that
    \begin{equation}\label{con12}
    H(1,2,1)+H(3,1)\equiv -\frac{4}{5}pB_{p-5}\pmod{p^2},
   \end{equation}
which substituting into (\ref{con8}) implies the congruence (\ref{con3}).\\
     
{\it Proof of {\rm (\ref{con4})}}.
Since $H_k=H_{k-1}+1/k$, for each $k=1,2,\ldots ,p-1$ 
the binomial formula implies
    \begin{equation*}\begin{split}
H_k^4-H_{k-1}^4
&=4\frac{H_{k-1}^3}{k}+
6\frac{\left(H_{k}-\frac{1}{k}\right)^2}{k^2}+
4\frac{H_{k-1}}{k^3}+\frac{1}{k^4}\\
&= 4\frac{H_{k-1}^3}{k}+6\left(\frac{H_k^2}{k^2}-
2\frac{H_{k-1}+\frac{1}{k}}{k^3}+\frac{1}{k^4}\right)+
4\frac{H_{k-1}}{k^3}+\frac{1}{k^4}\\
&=4\frac{H_{k-1}^3}{k}+6\frac{H_{k}^2}{k^2}-
8\frac{H_{k-1}}{k^3}-\frac{5}{k^4}.\nonumber
  \end{split}\end{equation*}
After summation of the above identity over $k$, 
and using the fact that by Wolstenholme theorem, 
$p^2\mid H(1)$,  we get
  \begin{equation}\label{con13}
4\sum_{k=1}^{p-1}\frac{H_{k-1}^3}{k}+
6\sum_{k=1}^{p-1}\frac{H_k^2}{k^2}-
8\sum_{k=1}^{p-1}\frac{H_{k-1}}{k^3}-5\sum_{k=1}^{p-1}\frac{1}{k^4}=
H_{p-1}^4\equiv 0\pmod{p^2}.
   \end{equation}
Similarly, we have
 \begin{equation*}\begin{split}
\frac{H_{k-1}^3}{k}&=\frac{1}{k}
\left(H_{k}-\frac{1}{k}\right)^3=\frac{H_k^3}{k}-3\frac{H_k^2}{k^2}+
3\frac{H_{k-1}+\frac{1}{k}}{k^3}-\frac{1}{k^4}\\
&=\frac{H_k^3}{k}-3\frac{H_k^2}{k^2}+
3\frac{H_{k-1}}{k^3}+\frac{2}{k^4},
  \end{split}\end{equation*}
which gives 
 \begin{equation}\label{con14}
\sum_{k=1}^{p-1}\frac{H_{k-1}^3}{k}=\sum_{k=1}^{p-1}\frac{H_k^3}{k}-3
\sum_{k=1}^{p-1}\frac{H_k^2}{k^2}+3\sum_{k=1}^{p-1}\frac{H_{k-1}}{k^3}
+2\sum_{k=1}^{p-1}\frac{1}{k^4}.
  \end{equation}
Substituting the expression (\ref{con14})   
into the first term of (\ref{con13}), we immediately obtain
  \begin{equation}\label{con15}
4\sum_{k=1}^{p-1}\frac{H_k^3}{k}-6\sum_{k=1}^{p-1}\frac{H_k^2}{k^2}+
4\sum_{k=1}^{p-1}\frac{H_{k-1}}{k^3}+3\sum_{k=1}^{p-1}\frac{1}{k^4}
\equiv 0\pmod{p^2}.
   \end{equation}
Further, since by  \cite[Theorem 3.2, p. 88]{z} 
$H(1,3)\equiv -\frac{9}{10}pB_{p-5}\,(\bmod{\,p^2})$, we have   
  \begin{equation}\label{con16}
\sum_{k=1}^{p-1}\frac{H_{k-1}}{k^3}=
\sum_{1\le i<k\le p-1}\frac{1}{ik^3}=H(1,3)\equiv 
-\frac{9}{10}pB_{p-5}\pmod{p^2}.
  \end{equation}
Inserting the  congruence (\ref{con11}) and the well known congruence 
$H(2)\equiv 0\,(\bmod{\,p})$  (see e.g., \cite[Theorem 3]{b})
into the shuffle relation $H(2)H(2)=2H(2,2)+H(4)$, we find that
 \begin{equation}\label{con17}
\sum_{k=1}^{p-1}\frac{1}{k^4}=H(4)\equiv 
-2H(2,2)\equiv \frac{4}{5}pB_{p-5}\pmod{p^2}.
      \end{equation}
Substituting (\ref{con3}), (\ref{con16}) and (\ref{con17}) into 
(\ref{con15}), we obtain the congruence (\ref{con4}).\\

{\it Proof of {\rm(\ref{con5})}}. The congruences (\ref{con16}) and  (\ref{con17})
yield
  $$
\sum_{k=1}^{p-1}\frac{H_k}{k^3}=
\sum_{k=1}^{p-1}\frac{H_{k-1}}{k^3}+
\sum_{k=1}^{p-1}\frac{1}{k^4}=H(1,3)+H(4)\equiv 
-\frac{1}{10}pB_{p-5}\pmod{p^2},
  $$
as desired.
This completes the proof of Theorem 1.1.
\end{proof}

\end{document}